\documentclass[11pt]{article}
\usepackage{amscd,amsthm,amssymb,amsfonts,amsmath,epsfig}
\theoremstyle{plain}
\newtheorem{thm}{Theorem}[section]
\newtheorem{lemma}[thm]{Lemma}
\newtheorem{prop}[thm]{Proposition}
\newtheorem{cor}[thm]{Corollary}
\theoremstyle{definition}

\theoremstyle{remark}
\newtheorem*{remark}{Remark}
\newtheorem*{ack}{Acknowledgments}

\newcommand{\Z}{\mathbb Z}    % Integers
\newcommand{\R}{\mathbb R}    % Real
\newcommand{\C}{\mathbb C}    % Complex
\newcommand{\PP}{\mathbb P}   % For projective space
\newcommand{\<}{\langle}   %\< is not defined yet. 
\renewcommand{\>}{\rangle} %\> is already defined.
\newcommand{\dT}{\partial T}
\newcommand{\dD}{\partial\Delta}
\newcommand{\DD}{\Delta^\vee_\gamma}
\newcommand{\dDD}{\partial\Delta^\vee_\gamma}
\newcommand{\orr}{\{ 0 \} }
\newcommand{\W}{\widetilde{W_\tau}}
\newcommand{\U}{\widetilde{U_\tau}}
\newcommand{\X}{X_\Delta}
\newcommand{\XX}{X_{\Delta_\gamma^\vee}}

\begin{document} 
\title{Torus Fibrations of Calabi-Yau
Hypersurfaces in Toric Varieties and Mirror Symmetry} 
\author{Ilia Zharkov\\Department of Mathematics\\ University of
  Pennsylvania\\Philadelphia, PA 19104\\izharkov@math.upenn.edu} 
\maketitle
  \begin{abstract}
  We consider regular Calabi-Yau hypersurfaces in $N$-dimensional
   smooth toric varieties. On such a hypersurface in the neighborhood of the
   large complex structure limit point we construct a fibration over a
   sphere $S^{N-1}$ whose generic fibers are tori $T^{N-1}$.  Also for
   certain one-parameter 
   families of such hypersurfaces we show that the monodromy
   transformation is induced by a translation of the $T^{N-1}$ 
   fibration by a section. Finally we construct a dual fibration
   and provide some evidence that it describes the mirror family.
  \end{abstract}

\section{Introduction} 
 
 Strominger, Yau and Zaslow [SYZ]  conjectured that
any Calabi-Yau manifold $X$ having a mirror partner $X^\vee$ should
admit a special 
Lagrangian fibration $ \pi: X \rightarrow B$ (a mathematical account
of their construction can be found in [M]). If so, the mirror
manifold $X^\vee$ is obtained by finding some suitable compactification of the
moduli space of flat $U(1)$-bundles along the nonsingular 
fibers, which restricts the fibers to be tori. More precisely, if
$B_0\subseteq B$ is the largest set such 
that $\pi_0=\pi\arrowvert_{\pi^{-1}(B_0)}$ is smooth, then $X^\vee$
should be a compactification of the dual fibration $R^1\pi_{0*}(\R/\Z)
\rightarrow B_0$.

The conjecture is trivial in the elliptic curve case. On a K3
surface the hyperk\"ahler structure translates the theory of special
Lagrangian $T^2$-fibrations to the standard theory of elliptic
fibrations in another complex structure. However, very
little progress has been made in higher dimensions so far, though
Gross and Wilson have 
worked out some aspects of the conjecture for the Voisin-Borcea
3-folds of the form $(K3\times T^2)/\Z_2$ [GW]. But the general question of
finding special Lagrangian fibrations on Calabi-Yau's still remains
open. 

We restrict our attention to the case of regular anticanonical hypersurfaces
in smooth toric varieties. The main result of our this paper is that
such a hypersurface in a neighborhood of the large complex structure
admits a torus fibration over a sphere. Unfortunately, we were unable 
to control the fibers to be special Lagrangian. However we will
argue that on some open patches the fibers do possess some calibration
property.

Batyrev [B] showed that toric varieties $X_{\Delta}$ with
ample anticanonical bundles are given by
reflexive polyhedra. Such a polyhedron $\Delta$ contains a unique integral
interior point \{0\}. A Calabi-Yau hypersurface $Y\subset X_{\Delta}$
is defined by an equation in the form $\sum_{\omega\in \Delta(\Z)}a_\omega
x^\omega =0$, where $\omega$ runs over the integral points in
$\Delta$. The image of $Y$ under the moment map $\mu : X_\Delta
\rightarrow \Delta$ has the shape of an amoeba (cf. [GKZ], Ch.6), a blob with
holes around some lattice points $\omega$ in $\Delta(\Z)$. The sizes of
the holes are determined by the 
corresponding coefficients $a_\omega$. If $Y$ is near the large complex
structure, $a_{\{0\}}$ is large, and $\mu(Y)$ has exactly one interior
hole corresponding to $\orr$, that is $\mu(Y)$ is homeomorphic to
$S^{N-1} \times I$. The idea is to choose 
the right trivialization of this product, so that for general $s\in S^{N-1}
\simeq \dD$ the preimage of the interval $\mu^{-1}(\{s\}\times I)$
would be an $(N-1)$-dimensional torus $T^{N-1}$. Singular fibers come
from the intervals $I_s:= \{s\} \times I$ which intersect the
$(N-2)$-dimensional skeleton of $\Delta$. For example, consider $Y$, a $ K3$
surface given by 
a quartic in $\C P^3 = X_\Delta$, where $\Delta$ is the integral
3-simplex in $\Z ^3$ with the vertices (-1,-1,-1), (-1,-1,3),
(-1,3,-1), (3,-1,-1). There
are exactly 4 points in $\mu(Y)$ on each of the 6 one-dimensional edges
of $\Delta$, corresponding to 4 points of intersection of $Y$ with the
projective line determined by this edge. Altogether they give 24
singular fibers. 

But instead of pursuing this idea we will modify the moment map
and  deform the original hypersurface. 
An explicit parameterization of the fibers will allow
us to analyze the action of the monodromy for one-parameter families of
hypersurfaces and construct a dual fibration. We will speculate that
this dual fibration represents the mirror Calabi-Yau.

Let us demonstrate almost all essential ideas by a simple example. Consider
a family of elliptic curves $E_t$ in $\C\PP^2$ given by the equations:
$$ txyz+x^3+y^3+z^3=0,$$
where $t$ plays the r\^ole of a parameter. As $t\rightarrow\infty$ the
curve $E_t$ degenerates to 3 lines with normal crossings. The
main idea is roughly to consider the asymptotic
behavior of $E_t$ up to the next order to keep the curve smooth. 

There are 6 regions in $E_t$
according to its image under the moment map (see Fig.1). In each of them
there are different terms in addition to $txyz$ which are dominant. 
For instance, in $U_{z}$ the elliptic curve $E_t$ for large $t$ is
approximated by 
$txyz+z^3=0$, and in $U_{zx}$ by $txyz+x^3+z^3=0$, etc. It is easy to
introduce a coordinate on a curve (which is still smooth in the
corresponding region) defined by the abbreviated equation. In $U_z$
either $x/z$ or $y/z$ is a coordinate, in
$U_{zx}$ we can use $x/z$ or $z/x$ and similar in the other
regions. The circle fibration is provided by fixing an absolute value
of the coordinate. The set of all possible
absolute values in all 6 regions clearly forms a circle, the base of
the fibration. The partition of unity technique of gluing these 6
pieces into one curve which approximates the original
elliptic curve constitutes section 3. 

\begin{center}
  \epsfig{file=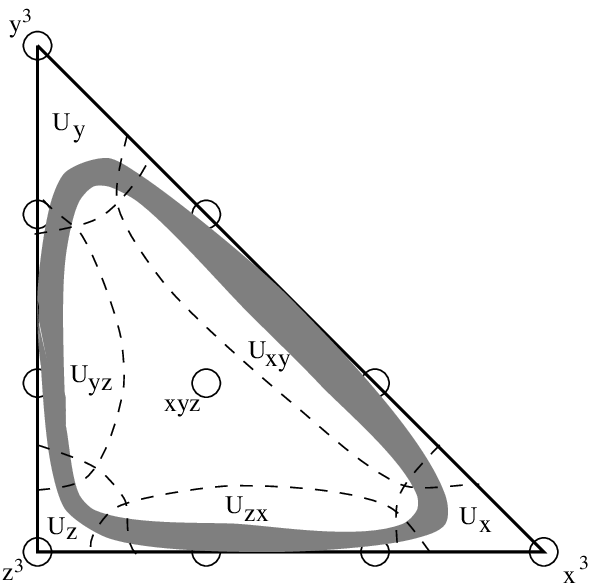}

{\scriptsize {\bf Fig.1} The image of the elliptic curve $
  txyz+x^3+y^3+z^3=0$ (shaded area) in $\C\PP^2$ under the moment map.}
\end{center}

To compute the monodromy as $\arg(t)\mapsto \arg(t)+2\pi i$, we need to
understand how the identification of the fibres changes in the
overlaps. Nothing happens in $U_{zx}$, $U_{xy}$ or $U_{yz}$, as the 
parameterization of fibers changes by reversing an orientation of the
circles, which does not depend on $t$ at all. Monodromy is nontrivial
only in  $U_{z}$, $U_{y}$ or $U_{x}$. Consider, e.g. $U_z$, where two
coordinates can be used $y/z$ or $x/z$, which are related by the
equation $t(x/z)(y/z)=-1$. As $\arg(t)\mapsto \arg(t)+2\pi i$, the
circles parameterized by $\arg(x/z)$ and by $\arg(y/z)$ are twisted by
$2\pi i$ with respect to each other. Combining all together we get a
triple Dehn twist. More careful monodromy calculations are performed
in section 4. 

\begin{center}
  \epsfig{file=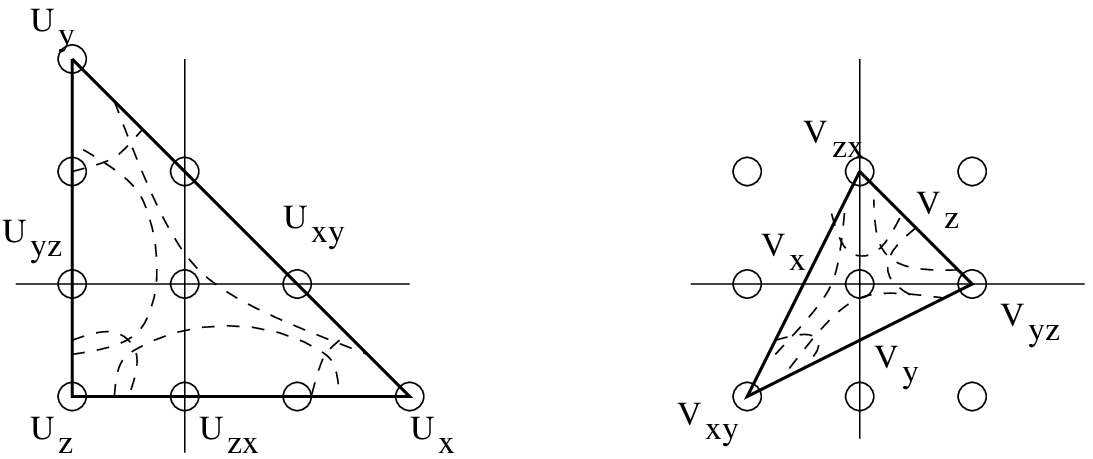}

{\scriptsize {\bf Fig.2} The mirror pair of elliptic curves.}
\end{center}

Section 5 is devoted to dual fibrations and mirror symmetry.
Leung and Vafa [LV] describe the idea of the mirror construction as
follows. $T$-duality interchanges small circle fibers in the corner regions
of one family with large circle fibers in the facet regions of the
mirror family. More precisely, let us consider the polytope
$\Delta^\vee$ dual to the above family of elliptic curves. In 
this case it is just the polar polytope $\Delta^D$ ,
hence again reflexive. The associated toric variety $X_{\Delta^D}$ is
singular but the elliptic curve in $X_{\Delta^D}$ is smooth because
it misses all singular points (the vertices of the triangle).

$\Delta^D$ is combinatorially dual to $\Delta$ and the dual elliptic
curve also breaks into 6 regions according to its image under the
moment map (see Fig.2). We will speculate that this dual curve also admits
similar fibration by circles and the fibers in the corresponding
regions $U_\alpha$ and $V_\alpha$ are naturally dual circles.

This example captures all the features except for the fact that in
higher dimensions there will always be singular fibers.

\begin{ack}
  I am very grateful to Ron Donagi, my thesis advisor, for suggesting this
  problem to me and constant supervising my work on it, and to Tony
  Pantev for illuminating discussions.
\end{ack}

\section{Hypersurfaces in toric varieties} 
 
 In this section we review some basic constructions in the theory of
 toric varieties and set up the notations. For more details see,
 e.g. [Cx] and references provided there.

Let $M\simeq\R ^{N}$ be an $N$-dimensional real affine space and $\Z
^{N}$ an integral lattice in $M$. We will choose an integral point
$\orr$ in $\Z^N$ to be the origin. This endows $M$ with a vector space
structure and  $\Z^N$ becomes a free abelian group.

An $N$-dimensional convex integral polyhedron $\Delta$ in $M$, is called {\it
 reflexive} if it contains the origin $\orr$ as an interior point and
 if its polar polyhedron
$$\Delta^D=\{u\in M^*:\<u,m-\orr\>\geq-1 \text { for all }
m\in\Delta\}\subset M^*$$ is also integral. We will denote by
$\Delta(\Z)$ the lattice points in the polyhedron $\Delta$, and by $\dD$
its boundary. It follows that each
$(N-1)$-dimensional face $\Sigma$ (which we will call a {\it facet} in
the future) of $\Delta$ is defined by an equation $\<u_\Sigma,m\>=-1$
for some $u_\Sigma\in M^*$. This easily implies that $\orr$ is the
unique integer interior point. For an arbitrary face $\Theta \subseteq
\dD$ we will denote by $\Theta_\Z$ the affine sublattice of $\Z^N$
generated by the integer points of $\Theta$.
    
Also we assume that $\Delta$ is {\sl nonsingular}. This mean that
every vertex of $\Delta$ is $N$-valent, that is exactly $N$ edges
$e_1,... , e_N$ emanate from it, and the integer points of these edges
${(e_i)}_\Z$, $i=1,...,N$, generate the lattice $\Z^N$. Because $\Delta$
is reflexive, i.e. the integral distance 
from the origin $\orr$ to any facet $\Sigma\subset\Delta
$ is 1, the lattice $\Z^n$ is as well generated by $\Sigma_\Z$
together with $\orr$. 

Let us denote by $X_\Delta$ the projective toric variety corresponding to
$\Delta$. The normal projective embedding is given, e.g., by the
closure of the image of the map $(\C^*)^N \hookrightarrow
\PP^{|\Delta(\Z)|-1}$, $x\mapsto\{x^{\omega_1}:x^{\omega_2}:
... :x^{\omega_{|\Delta(\Z)|}}\}$. Because $\Delta$ is nonsingular, the
toric variety $\X$ is smooth.
We will use $\{x^{\omega}\}$ as projective coordinates on
$\X$. One may want to
restrict the set of monomials $\{x^{\omega}\}$ to a subsystem of the
anticanonical linear system. In this case we allow $\omega$ vary among  
$A\subset\Delta(\Z)$, a subset of integer points in $\Delta$, such that
$\orr\in A$ and $\Delta$ is the convex hull of $A$.

The moment map $X_\Delta \rightarrow \Delta$ is given
by 
$$ \mu(x)=\frac {\sum_{\omega\in A} |x^\omega|\cdot
\omega}{\sum_{\omega\in A} |x^\omega|}.$$ 
It is a well-defined function on $X_\Delta$, because both the top and bottom
are polynomials of the same homogeneous degree. 

A {\it triangulation} of a convex polyhedron is a decomposition of it
into a finite number of simplices such that the intersection of any
two of these simplices is a common face of them both (maybe empty). By
a triangulation $T$ of $(\Delta,A)$ we simply mean a triangulation of
$\Delta$ with vertices in $A$. Note that we do not require every
element of $A$ to appear as a vertex of a simplex. A continuous
function $\psi:\Delta\rightarrow\R$ is called {\it T-piece-wise
  linear} if it is affine-linear on every simplex of $T$. Such a
function $\psi$ is {\it convex} if for any $x,y\in\Delta$, we have
$\psi(tx+(1-t)y)\geq t\psi(x)+(1-t)\psi(y), 0\leq t\leq 1$. We call it
{\it strictly convex} if the (maximal-dimensional) simplices of $T$
are the maximal domains of linearity of $\psi$.
 
A triangulation $T$ of $(\Delta,A)$ is called {\it coherent} (some
 authors call it regular or projective) if there exists a strictly convex
 $T$-piece-wise linear function. Call a coherent triangulation $T$
 {\it central} if $\orr$ is a vertex in every ( maximal dimensional)
 simplex of $T$. Let $\dT$ denote the collection of  
simplices of $T$ lying in $\dD$. We will use $\sigma$ to denote
maximal dimensional simplices in $\dT$ and $\tau$ for arbitrary simplices
of $\dT$. Denote by $C_\tau\subset\Delta$ the corresponding simplex
 in $T$ of one dimension higher with the base $\tau$ and the vertex $\orr$.
We will denote by $\tau_\Z$ the sublattice of $\Z^N$
 generated by the 
 integer points in $\tau$. $\Lambda_\tau$ will denote the sublattice in
 $\tau_\Z$ of index $(\dim\tau)!\cdot vol(\tau)$ generated by the
 vertices of $\tau$.

Given a triangulation $T$ of $(\Delta,A)$, every function
  $\lambda_\R:A\rightarrow\R$ defines a {\it characteristic} function
  $\psi_{\lambda_\R}:\Delta\rightarrow\R$, a $T$-piece-wise linear
  function, by its values on the vertices of $T$. Denote by
  $C(T)\subset\R^{|A|}$ a subset of such functions $\lambda_\R$, whose
  corresponding characteristic functions $\psi_{\lambda_\R}$ are
  convex and  $\psi_{\lambda_\R}(\omega)\geq \lambda_\R(\omega)$ for
  any $\omega\in A$. A person familiar with toric varieties
  immediately recognizes a secondary cone, the normal cone to the
  secondary polyhedron at the vertex corresponding to the
  triangulation $T$. In particular $C(T)$ has non-empty interior if
  $T$ is a coherent triangulation. The last piece of data we will need
  is an integral vector $\lambda$ in the interior of $C_T \subset
\R^{|A|}$. This means that the characteristic function $\psi_\lambda: \Delta
\rightarrow \Z$ is strictly convex with respect to the triangulation
  $T$ and  $\psi_\lambda(\omega)\geq \lambda(\omega)$, with equality
holding exactly for the vertices of $T$. 

From now on we fix the following data: a nonsingular reflexive
integral polyhedron $\Delta$, a subset of its integral points $A$, a
central coherent 
triangulation $T$ and an integral vector $\lambda\in C(T)$. Given such
data we can define the 1-parameter family of Calabi-Yau hypersurfaces
$F_t$  by the equations
$$t^{\lambda(0)}x^{\orr}-\sum_{\omega\in A\cap\dD}
t^{\lambda(\omega)}x^\omega=0.$$ 

To any hypersurface given by an equation in the form  $\sum_{\omega\in
 A}a_\omega x^\omega =0$ we can associate the vector
 $\lambda_a:=\{\log|a_{\omega_1}|,...,\log|a_{\omega_{|A|}}|\}\in\R^{|A|}$. 
 If this vector $\lambda_a$ is sufficiently far from the walls of the
 secondary fan then the corresponding hypersurface in $\X$ is
 nonsingular. Note that for $|t|\gg 0$,
$\log | t^{\lambda(\omega)}| = \lambda \cdot
\log|t|$ lies deeply inside $C_T$, so that the hypersurface $F_t$ is
smooth (cf. [GKZ], Ch 10). The main object of study will be the
behavior of this family as $|t|\rightarrow\infty$. In the limit we
 approach  the large complex structure limit point.

It may sometimes be convenient to have local coordinates on
affine subsets of $\X$. We define $\tau-associated$ coordinates for any
$k$-dimensional simplex $\tau\in\dT$ in the triangulation. Choosing
coordinates on the open orbit $(\C^*)^N\subset\X$ is equivalent to
choosing an affine basis 
$\{\Omega,\omega_1,...,\omega_N\}$ for the lattice $\Z^N$ together
with a marked reference point
$\Omega$.  Consider the minimal face $\Theta_\tau\subset\Delta$
containing $\tau$. Let $l$ be its  dimension, $k\leq l\leq N-1$.
We use the fact that $\Delta$ is reflexive and nonsingular to choose
$\{\Omega,\omega_1,...,\omega_N\}$ such that the reference point
$\Omega$ is a vertex of $\tau$, the first $k+1$ vertices
$\{\Omega,\omega_1,...,\omega_k\}$ generate the lattice $\tau_\Z$, the
first $l+1$ vertices $\{\Omega,\omega_1,...,\omega_l\}$ generate the
lattice $(\Theta_\tau)_\Z$. Moreover, because of the unit integral distance
from the facets of $\Delta$ intersecting at $\Theta_\tau$ to the origin
$\orr$, we can have the rest of the basis satisfy
$\orr-\Omega=(\omega_{l+1}-\Omega)+...+(\omega_N-\Omega)$. Notice that
the last $N-l$ points are uniquely determined modulo $(\Theta_\tau)_\Z$. Then
$\{y_1,...,y_N\}:=\{x^{\omega_1-\Omega},...,x^{\omega_N-\Omega}\}$
give local coordinates on $(\C^*)^N$. In fact, they extend to
coordinates on the affine subset of $X_\Delta$ obtained by removing the
divisors which correspond to the facets of $\Delta$ not containing
$\Omega$. 

It may also be useful to identify $\Delta$ with the closure of the
positive real part ${(\X)}_{\geq 0}$ of $\X$. The homeomorphism ${(\X)}_{\geq
0}\simeq\Delta$ is provided by the restriction of the moment map. In
particular, $\{|y_1|,...,|y_N|\}$ will provide coordinates in $\Delta$
with facets not containing $\Omega$ excluded as above.

Inspired by the proof of Viro's theorem (see [GKZ], Ch. 11) we will
use the weighted moment map $\mu_t: X_\Delta \rightarrow \Delta $ defined as
$$
\mu_t(x)=\frac {\sum_{\omega\in A} |t^{\lambda(\omega)}|
\cdot|x^\omega|\cdot \omega}{\sum_{\omega\in A}
|t^{\lambda(\omega)}| \cdot|x^\omega|}.
$$

The right way to think about this weighted moment map is the
following. Add one extra dimension to $M \cong \R^N$ and
extend the lattice to $\Z^{N+1} \subset \R^{N+1}$. Let $P$ be the
convex hull of $\{(\omega,\lambda(\omega))\}_{\omega \in A}$ in
$\R^{N+1}$ (see fig.3).  Then
we can think of the whole family $\{F_t\}$ as a hypersurface in $X_P$ ($t$
is considered as a coordinate). The vertical projection  $p : P
\rightarrow \Delta$ splits the boundary of $P$  into two pieces
 $\partial P=\partial_+ P \cup \partial_- P$. In fact $\partial_+P$ is
 exactly the graph of the characteristic function $\psi_\lambda$ and
 the projection $p:\partial_+ P\rightarrow\Delta$ identifies the faces
 of $\partial_+ P$ with the simplices in the triangulation $T$. The
 weighted moment map $\mu_{t_0}$ will be just the composition of the
 restriction $\mu^{(P)}_{t_0} := \mu^{(P)} \arrowvert_{t=t_0}:
 {X_P}\arrowvert_{t=t_0} \rightarrow P$ of the ordinary moment map
 $\mu^{(P)}$ to the hypersurface $H_{t_0}:=\{t=t_0\}$ in $X_P$ with
 the vertical projection $p : P \rightarrow \Delta$.

The image of $\mu^{(P)}(H_{t_0}) \subset P$ is the graph of a function
$\Psi_{t_0}: \Delta \rightarrow \R$ with the following
crucial property (see [GKZ], Ch. 11):

\begin{prop}
  As $|t_0|\rightarrow \infty$ the function $\Psi_{t_0}: \Delta
  \rightarrow \R$ is continuous and smooth outside $\dD$, converging
  uniformly to the characteristic function $\psi_\lambda$, whose graph
  is $\partial_+ P$.
\end{prop}

\begin{center}
  \epsfig{file=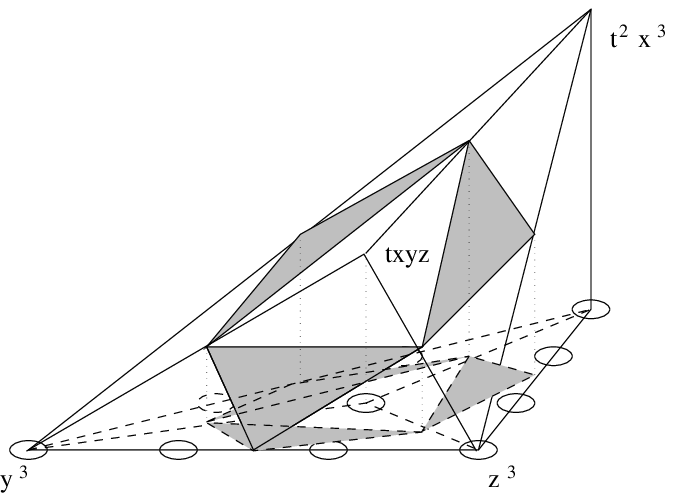}

     {\scriptsize  {\bf Fig.3} The extended polyhedron $P$ for the
     family $txyz-t^2x^3-y^3-z^3=0$ in $\C\PP^2$ and the image of the
     weighted moment map as $t\rightarrow\infty$.}
\end{center}

The image of the hypersurface $\{F_t\}\subset X_P$ under the moment
map $\mu^{(P)}$ 
misses all the vertices of $P$ (cf. [GKZ], Ch.6) and, in particular, (because
$\psi_\lambda$ is convex) it misses some neighborhood of $(\orr,
\lambda(0)) \in P$. Applying the above proposition we see that for
sufficiently large $|t|$ the image of the weighted moment map
$\mu_t(F_t) \subset \Delta$ misses some ball $B$ around $\orr \in
\Delta$. Denote by $\Delta^0 :=\Delta \backslash B$ our polytope with
that ball removed. Clearly $\Delta^0$ is homeomorphic to
$S^{N-1}\times I$ and the next step will be to find a good
trivialization of this product.

Let us conclude this section with a list of notations
used throughout the rest of the paper.

$M=\Z^N\otimes\R\simeq\R^N$ the real affine space;

$\Z^N$ an $N$-dimensional lattice in $M$.

$\Delta\subset M$ a convex nonsingular integral reflexive polyhedron,
$\dD$ its boundary;

$\Delta(\Z)$ the set of integral points in $\Delta$;

$\orr\in\Delta(\Z)$ the unique integral interior point;

$A$ a subset of $\Delta(\Z)$ containing $\orr$ and all vertices of
$\Delta$;

$\Theta$ a face of $\Delta$, $\Sigma$ a maximal dimensional face of
$\Delta$;

$T$ a central coherent triangulation of $(\Delta,A)$;

$\lambda\in C(T)$ an integral vector in the interior of the secondary
cone at $T$;

$\dT$ the induced triangulation of $\dD$;

$\tau$ a $k$-dimensional simplex in $\dT$, $\sigma$ a maximal
dimensional simplex in $\dT$;

$\Theta_\tau$ the minimal face of $\Delta$ containing $\tau$;

$C_\tau$ the $(k+1)$-dimensional simplex in $T$ over $\tau$ with the
vertex $\orr$;
 
$O(\tau)$ the center of a simplex $\tau$;

$\Theta_\Z$ or $\tau_\Z$ the affine sublattices of $\Z^N$ generated by
the integral points in $\Theta$ or $\tau$;

$\Lambda_\tau$ the sublattice of $\tau_\Z$ generated by the vertices
of $\tau$;

$\X$ the toric variety associated to $\Delta$;

$F_t$ the family of the Calabi-Yau hypersurfaces in $\X$; 

$\mu_t$ the weighted moment map;

$\Delta^0$ the polyhedron $\Delta$ with a small ball around $\orr$ removed,
$\X^0:=\mu_t^{-1}(\Delta^0)$.

\section{Torus fibrations}

At this moment unfortunately we must leave the realm of beautiful
algebraic geometry and employ some analysis techniques like partitions
of unity and transversality theory.  Let $Bar(\dT)$ be the first
barycentric subdivision of the 
triangulation $\dT$. The vertices in this subdivision are the centers
$O(\tau)$ of the simplices $\tau$ in $\dT$. Consider the subdivision
$\{V^0_\tau\}_{\tau\subset\partial T}$ {\it dual} to $Bar(\dT)$. Namely,
take the second barycentric subdivision $Bar^{(2)}(\dT)$ of $\dT$ and
define $V^0_\tau$ to be the union of all simplices in $Bar^{(2)}(\dT)$
having $O(\tau)$ as a vertex. Every $V^0_\tau$ contains the point
$O(\tau)$ for a unique $\tau$ and is
labeled correspondingly (see Fig.4). For each $V^0_\tau$ we take its
small open neighborhood 
$V_\tau$ to get an open cover of $\dD$. By construction $V_\tau$ and
$V_{\tau'}$ intersect iff either $\tau\subset\tau'$ or
$\tau'\subset\tau$. So every point in $\dD$ lies in at most $N$
different $V_{\tau_i}$'s for $\{\tau_i\}$ forming a nested sequence. 

\begin{center}
   \epsfig{file=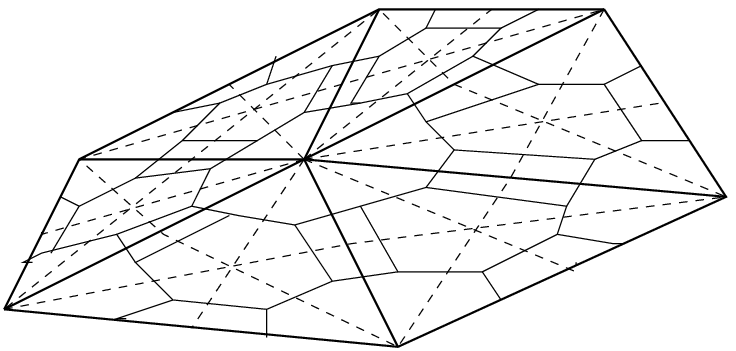}

     {\scriptsize {\bf Fig.4} $V^0_\tau$ subdivision of $\dD$.}
\end{center}

Define the following subsets of $\dD$:
$$
U_\tau:=V_\tau - \bigcup_{\tau'\neq\tau}\overline V_{\tau'}, \quad
W_\tau:=\overline V_\tau - \bigcup_{\tau'\supset\tau}V_{\tau'},\quad \text
{(set-theoretic difference)},
$$ 
where $\overline V_{\tau}\subset\dD$ denotes the closure of $V_\tau$.

\begin{center}
  \epsfig{file=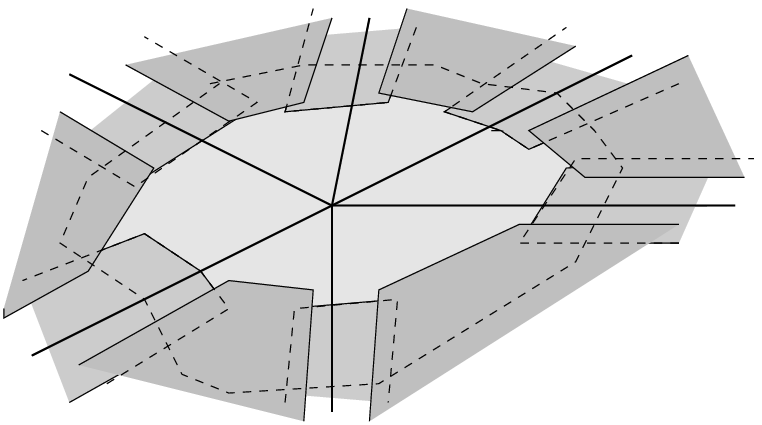}

     {\scriptsize  {\bf Fig.5} $W_\tau$-cell decomposition of $\dD$.}
\end{center}

The collection of $cells$ $\{W_\tau\}_{\tau\in\dT}$ provides a
CW-decomposition of $\dD$ homeomorphic to that given by $\{V^0_\tau\}$
 (see fig.5).  $U_\tau$ are the ``pure'' open subsets in
$W_\tau\subseteq \overline V_\tau$. 

Now we will construct a trivialization of $\Delta^0 \simeq \dD\times
I$. This will provide $\X^0:=\mu_t^{-1}(\Delta^0)$ with the structure of
a fibration via the composition of the maps
$$\X^0\stackrel{\mu_t}{\longrightarrow}\Delta^0\simeq\dD\times
I\stackrel{pr_1}{\longrightarrow}\dD.$$ 
Inside the small central ball $B^0$ choose
a concentric mini-copy $\Delta'$ of $\Delta$  with the induced 
triangulation $\dT'$ of $\dD'$. To construct a trivialization we need
to connect the two boundaries $\dD$ and $\dD'$ by non-intersecting
intervals. To do this we need a rule how to choose which pairs of points
$s\in\dD$ and $s'\in\dD'$ are to be connected and an interval
connecting them. After that we must make sure that these intervals do
provide a trivialization of $\Delta^0$. We will denote by
$\Theta'\subset\dD'$ and $\tau'\in\dT'$ the mini-copies of $\Theta$
and $\tau$ correspondingly. 

Let $s'$ be an interior point of an $l$-dimensional face
$\Theta'\subset\dD'$. We choose $\tau$-associated coordinates
$\{y_1,...,y_N\}$ for some $\tau$, such that $\Theta_\tau$ is the
minimal face containing $\tau$.  Then $\{|y_1|,...,|y_N|\}$ will
provide coordinates in the open subset of $\Delta$ corresponding
to $\tau$ (see the previous section) by means of the weighted moment
map. Let $m_i=|y_i|(m)$ be  
the coordinates of a point $m\in\Delta$. Define {\it the curved normal
cone} to $\Theta'$ at the point $s'\in\Theta'$ by
$$n(s')=\{m\in\Delta:m_i=s'_i,\ i=1,...,l,\text{ and } m_i\leq s'_i,\
i=l+1,...,N\}.$$
 Note that the definition does not depend on the
choice of $\tau$-coordinates. Combined all together the curved normal
cones form a {\it fat curved normal fan} to the polyhedron
$\Delta'$. Namely, to every $k$-dimensional face $\Theta'\subset\dD'$
we can associate an $(N-l)\times l$-dimensional fat curved normal cone
$N(\Theta')\simeq n(\Theta')\times \Theta'$, where $n(\Theta')$ is the
curved normal cone to $\Theta'$ at any point $s'\in\Theta'$. This fat
fan provides a fat cone decomposition of $\Delta$ with $\Delta'$
deleted.

\begin{center}
  \epsfig{file=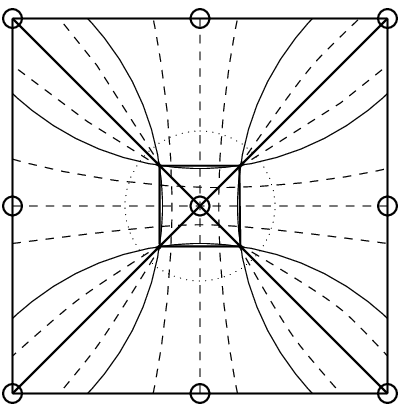}

{\scriptsize {\bf Fig.6} The fat cone decomposition and
  $I_s$-fibration of $\Delta^0$ for $\X=\C\PP^1\times \C\PP^1$.} 
\end{center}

The trivialization of $\Delta^0\simeq\dD\times I$ is achieved by
connecting any point $s'\in \dD'$ with all the points $s\in\dD$ lying
in the curved normal cone $n(s')$. The connecting intervals are given
by the unique line (given by affine linear equations in the
remaining $(N-l)$ coordinates) in
$n(s')$ passing through both points $s$ and $s'$ (see fig.6). Because
every point 
$s\in\dD$ belongs to a unique interval, we will denote this interval
by $I_s$. Notice that the interval $I_s$ does not depend on the choice of
$\tau$-coordinates, for a change of coordinates does not affect the
affineness of the defining equations in the last $N-l$ coordinates. 
But $I_s$ does depend
on the value of $t$ through the dependence of the weighted moment map
on $t$ and this $t$-dependence will be crucial in the next
proposition. 

For any $k$-dimensional simplex $\tau\in\dT$ we will refer to
$s\in\dD$ as a $\tau$-$point$ if in any $\tau$-associated coordinate
system the first $k$ coordinates $\{|y_1|,...,|y_k|\}$ are constant
along the interval $I_s$. In other words, $I_s$ lies in a fat cone
$N(\Theta')$ for some $\Theta'\supset\Theta'_\tau$. Note again that
this does not depend on the choice 
of the $\tau$-associated coordinates. We will say that a collection
$\{I_s\}_{s\in \dD}$ is $W$-$supported$ if for any $\tau\in\dT$ all
the points in the cell $W_\tau$ are $\tau$-points.

\begin{prop}
  Given a cell decomposition $\{W_\tau\}$ there exists a positive real
  number $R$, 
  such that for $|t|\geq R$  the collection $\{I_s\}_{s\in \dD}$ 
   providing a trivialization of $\Delta^0 \simeq
  S^{N-1}\times I$ is $W$-supported.
    \begin{proof}
    We claim that in the limit $t\rightarrow\infty$ the intervals
    $I_s$ will provide a bijective correspondence between points
    in $\tau$ and points in $\tau'$ for all $\tau\in\dT$. In particular,
    this means that all interior points in any  simplex $\tau$ eventually
    become $\tau$-points. Assuming this claim the
    proposition follows immediately from the following observation. $W_\tau$
    is compact and is contained in the interior of the union of the
    facets of the polytope $\Delta$ containing the face
    $\Theta_\tau$, and every point in this interior becomes a
    $\tau$-point for large $|t|$.

    To show the claim we consider the extended  polytope $P$. Keeping
    the first $k$ $\tau$-associated coordinates fixed is
    equivalent to fixing $(k+1)$ projective coordinates
    $\{|t^{\lambda(\omega)}x^\omega|\}_{\omega\in\tau}$. The 
    restriction of these equations to $p^{-1}(C_\tau)\subset\partial_+
    P$, where $p:\partial_+P\rightarrow\Delta$ is the vertical
    projection, defines a ray $R_s$ from the origin $(\orr, \lambda(0))$ to
    a point in $p^{-1}(\tau)$. Thus we see that every point in the interior
    of $p^{-1}(\tau')$ connects to a unique point in the interior of
    $p^{-1}(\tau)$. But because the function $\Psi_{t_0}: \Delta
    \rightarrow \R$ converges uniformly to the characteristic function
    $\psi_\lambda$ as $|t|\rightarrow\infty$, $I_s$ 
    would be given in the limit by the projection of the ray $R_s$ to
    $\Delta$. 
  \end{proof}
\end{prop}

\begin{remark}
  If one wants a fibration of $\X^0$ over a {\sl smooth} sphere
 as opposed to just topological, then some care should be taken to
 smooth out $I_s$ near the boundaries of the fat cone decomposition of
 $\Delta^0$. From now on we will assume that $\dD\simeq S^{N-1}$ is
 endowed with a smooth structure compatible with the fibration.
\end{remark}

For any subset $S\subset\dD$ we will denote by $I_S\simeq S\times I$
the union of $I_s$ for $s\in S$, and let $\widetilde S
:=\mu_t^{-1}(I_S)\subset\X^0$. In particular, we will be using $\U$
and $\W$. 

Next we choose a partition of unity  $\{\rho^0_\tau\}$
subordinate to the cover $V_\tau$, and define
the cut-off functions $\rho_\omega: \dD\rightarrow[0,1]$ for $\omega\in
\dT$, a vertex of the triangulation, by
$\rho_\omega:=\sum_{\tau\ni\omega}\rho^0_\tau$. In particular, it is clear
that
$$\rho_\omega(s) =
\begin{cases}
  0,\quad\text {unless $s\in V_\tau$ for some $\tau$ containing
  $\omega$,} \\ 1,\quad\text{if $s\in V_\tau$ only for those $\tau$ which
  contain $\omega$.}\\
\end{cases} $$
To uniformize the notation let $\rho_{\orr} \equiv 1$ and $\rho_
\omega\equiv 0$ for those $\omega\in A$ which are not vertices in the
triangulation. Extend $\rho$'s to the entire $\Delta^0$ by setting
$\rho_\omega (I_s):=\rho_\omega(s)$ for the entire interval
$I_s$. We will denote by $\rho_\omega$ also the pull back of the
cut-off functions
to $X_\Delta^0:= \mu_t^{-1}(\Delta^0)$ via the moment map.

Now we are in position to define the auxiliary object $H_t$, a real
(non-analytic) codimension 2 submanifold in $X_\Delta$, which we will
still call a hypersurface. It is defined by the following
equation:
$$ t^{\lambda(0)}x^{\orr}-\sum_{\omega\in A\cap\dD}\rho_\omega
    t^{\lambda(\omega)}x^\omega =0.
$$

 The restriction of $H_t$ to $\U$ defines an open set (in complex topology)
 of an algebraic subvariety in $X_\Delta$ given by the equation 
$$ t^{\lambda(0)}x^{\orr}-\sum_{\omega\in \tau}
t^{\lambda(\omega)}x^\omega =0.
$$
 And the cut-off functions $\rho$ are
designed to connect these pieces together.

 We showed that $\mu_t(F_t)$ misses a ball around $\orr\in\Delta$. A
 similar argument applies to show that $\mu_t(H_t)$
 also lies in $\Delta^0$. Namely, every point of $H_t$ can be thought
 as a point in an algebraic hypersurface in $\X$ defined by the same
 equation as $H_t$ but with constant $\rho$'s. Every such algebraic
 hypersurface clearly misses a ball around $\orr\in\Delta$. Using the
 compactness of $H_t$ we get the claim.

The structure of the rest of this section is
 the following. The fibration $\X^0\rightarrow\dD$, when restricted to
 the auxiliary hypersurface, induces a torus fibration $h_t:
 H_t\rightarrow\dD$. Then we will show that the auxiliary  hypersurface 
 is, in fact, diffeomorphic to the original one. Moreover this
 diffeomorphism extends to a diffeomorphism between the entire
 families $\{F_t\}$ and $\{H_t\}$. This will provide a torus
 fibration of the original Calabi-Yau hypersurface.

 The first step is to show that the collection $\{I_s\}_{s\in\dD}$ defines
 a torus fibration $h_t:H_t\rightarrow\dD$. Of course, some of the fibers are
 degenerate and we will try to describe them as explicitly as
 possible. Let $\tau$ be a $k$-dimensional simplex 
 in $\dT$, and $\Theta_\tau\subset\Delta$ the minimal face containing
 $\tau$, let $l:=\dim \Theta_\tau$. Let  $\{y_1,...,y_l\}$ be the first
 $l$ of $\tau$-associated coordinates. A $\tau$-point $s$ defines an
 $l$-dimensional torus  $T^l\subset(\C^*)^l$ by setting $|y_i|=s_i$.
 Let $P(y_1,...,y_k)=0$ be the equation 
$$\sum_{\omega\in \tau}\rho_\omega t^{\lambda(\omega)}x^\omega =0,$$
written in the local coordinates.
 Denote by $D_{\tau,s}\subset T^l$ the zero locus of the polynomial
 $P(y_1,...,y_k)$. To get some idea what $D_{\tau,s}$ looks like we consider a
 $k$-dimensional torus $T^k\subset(\C^*)^k=\text{Spec}[z_1^{\pm
   1},...,z_k^{\pm 1}]$, defined by fixing $|z|$. Denote by
 $D^0_{\tau,|z|}$ the intersection of $T^k$ with a plane
 $\rho_0+\sum_{i=1}^k\rho_i z_i=0$. For generic $|z|$ and $\rho$,
 $D^0_{\tau,|z|}$ will be either empty or a $(k-1)$-torus. In
 exceptional cases $D^0_{\tau,|z|}$ can be a single (real) point. The
 relation between  $D_{\tau,s}$ and  $D^0_{\tau,|z|}$ can be described
 as follows. The substitution
 $z_i=t^{\lambda(\omega_i)-\lambda(\Omega)} x^{\omega_i-\Omega}$, where 
 $\{\Omega,\omega_1,..., \omega_k\}$ are the vertices of the simplex
 $\tau$, defines a covering map $\pi_\tau:T^k\rightarrow T^k$. The
 degree of $\pi_\tau$ is equal to the index of the sublattice
 $\Lambda_\tau$ inside the lattice $\tau_\Z$, which is given by
 $k!\cdot vol(\tau)$. Then $D_{\tau,s}\subset T^l$ is a pull back of
 $D^0_{\tau,|z|}\subset T^k$ under the composition of maps
$$T^l\stackrel{pr}{\longrightarrow} T^k\stackrel{\pi_\tau}
{\longrightarrow}T^k,$$ where
$pr:\{y_1,...,y_l\}\rightarrow\{y_1,...,y_k\}$ is a projection onto
the first $k$ coordinates.

\begin{prop}
  Let $s\in W_\tau$ as above be a point in the interior of an
  $L$-dimensional face of $\Delta$, $L\geq l$. Then the fiber
  $T_s:=\mu_t^{-1}(I_s)$ itself has a structure 
  of a fibration $p_s:T_s\rightarrow T^l$ with a generic fiber
  $T^{N-1-l}$ and fibers $T^{L-l}$ over the discriminant locus
  $D_{\tau,s}$. Thus, $T_s$ is homeomorphic to 
$$T^{L-l}\times\left((T^l\times T^{N-1-L})/\sim\right),\quad
where\  (d,t_1)\sim (d,t_2),\text{ if }d\in D_{\tau,s}.$$
 In particular, if the point $s$ is in the interior of an $(N-1)$
dimensional face, then $T_s$ is a smooth $(N-1)$-torus.
  \begin{proof}
    We choose $\tau$-associated coordinates
    $\{y_1,...,y_N\}$. In particular, they provide coordinates in $\W$. The
    equation of the auxiliary hypersurface restricted to $T_s$ becomes:  
 $$ y_{l+1}y_{l+2}...y_N=P_t(y_1,...,y_k),$$
    where $P_t$ is a polynomial (all $\rho$ are constants on $T_s$).
    Because $s$ is a $\tau$-point, $|y_i|$ are fixed and non zero for
    $i=1,...,l$. This gives a projection $p_s:T_s\rightarrow T^l$. A fiber
    of this projection is determined by fixing a point $Y:=\{y_i\}$,
    $i=1,...,l$, on the base. After that we are left with the equation
    $ y_{l+1}...y_N=P_t(Y)=const$. At this point we must remember
    that $I_s$ is given by a line in chosen coordinates, hence it
    intersects the hyperbola $|y_{l+1}|...|y_N|=|P_t(Y)|$ in exactly
    one point in $\Delta$. Thus it determines the remaining $|y_i|$ uniquely.

   We see that a generic fiber of the projection  $p_s:T_s\rightarrow
    T^l$ is an $(N-l-1)$-torus. The dimension drops to $(L-l)$, if
    $P_t(Y)=0$, i.e. exactly if the point $Y$ is in the discriminant locus
    $D_{\tau,s}$. 
  \end{proof}
\end{prop}
  
We  want to say some words about the discriminant locus $D(H_t)$ of the 
fibration $h_t:H_t\rightarrow\dD$. The above proposition says that
$D(H_t)$ consists of all points in $\partial\dD$, the
$(N-2)$-dimensional skeleton of $\dD$, which are in the image of the
moment map $\mu_t$. Thus $D(H_t)$ is homeomorphic to 
$\partial\dD$ with some neighborhoods of its vertices removed,
which has a homotopy type of $Sk_T^{(N-3)}$, the $(N-3)$-skeleton of the
subdivision dual to the triangulation of $\partial\dD\subset\dD$. Moreover,
with an appropriate choice of the $V$-subdivision ($W_\tau$ should
have small volume for all $\tau$ with $1\leq\dim\tau\leq N-2$) the
discriminant locus $D(H_t)$ 
will lie in an arbitrarily small neighborhood of $Sk_T^{(N-3)}$. We
will refer to this limit as the {\it right} $W$-decomposition limit.

To deform the auxiliary hypersurface $H_t$ and to use 
transversality theory we must make sure that it is smooth.
\begin{lemma}
  For a generic choice of $\{\rho\}$ the hypersurface $H_t$ is smooth. 
  \begin{proof}
    It is enough to show that $H_t$ is smooth in
    $\W:=\mu_t^{-1}(I_{W_\tau})$ for every $\tau\subset\dT$. We will
    use $\tau$-associated  coordinates $\{y_1,...,y_N\}$ in
    $\W$. Let $G_t$ be the defining equation of the auxiliary
    hypersurface in $\W$ (see the  previous proposition) 
$$ G_t:=y_{l+1}y_{l+2}...y_N-b_\Omega\rho_\Omega(|y|)-\sum_{i=1}^k
b_i\rho_i(|y|) y^{\alpha_i}=0,$$ 
where $t$-dependence is encoded into $b_i$'s. 

Denote by ${\mathcal R}$ the family of $\{\rho_\omega\}$ constructed
from the family of partitions of unity subordinate to $\{V_\tau\}$.
  Consider the map $G_t^\rho: \W\rightarrow\C$. The statement
that $G_t^\rho=0$ is smooth in this language translates as $G_t^\rho$
is transversally regular to $0\in\C$. We want to show that there are
enough functions in ${\mathcal R}$, so that a generic  $G_t^\rho$
is transversally regular to $0\in\C$. According to the restricted
transversality theory, it is enough to show that the map of the entire
family $\widetilde{G_t}:\W\times{\mathcal R}\rightarrow\C$ is
transversally regular to  $0\in\C$ (cf, e.g. [DNF]). For any point
$\tilde{x}=(x,\rho)\in\widetilde{G_t}^{-1}(0)$ we need to show that
the tangent space at
$\tilde{x}$ maps onto $\C$. Consider the restriction of
$\widetilde{G_t}$ to the slice in a small neighborhood of $\tilde{x}$,
given by $\rho=const$. The function $\widetilde{G_t}
\arrowvert_{\rho=const}$ becomes algebraic and it is a straightforward
calculation to show that the tangent space to that slice is
transversal to $0\in\C$. Just notice that all $y_i$, $i=1,...,k$ are
non zero and at least one of the $\rho_i$ is non zero too. Hence a
generic choice of $\rho$ will provide a smooth preimage of
$0\in\C$. This completes the proof. 
  \end{proof}
\end{lemma}

The next step is to find a small
deformation of $H_t$ and a diffeomorphism of
$X_\Delta^0:=\mu_t^{-1}(\Delta^0)$ inside $X_\Delta$, which transforms
the deformed equation for $H_t$ into the equation of a genuine
hypersurface $F_{(\Gamma\cdot t)}$ for some real number $\Gamma$. For
this we need a technical lemma.

 \begin{lemma}
   There exists a function $\chi(s,\gamma,\omega): \dD\times\R_{\geq
   0}\times A\rightarrow \R$, smooth with respect to
   $(s,\gamma)$ and affine linear with respect to $\omega$, and
   satisfying:

    $ \bullet\ \chi(s,0,\omega)\equiv 0$ and
    $\chi(s,\gamma,\orr)\equiv \gamma\cdot\lambda(0)$.

    $\bullet$ As ${\gamma\rightarrow\infty}$ the function
    $e^{\gamma\lambda(\omega)-\chi(s,\gamma,\omega)}$ converges
    (uniformly) to $\rho_\omega (s)$ for every $\omega\in A$.
    \begin{proof}
     First,  for every $\tau \subset\dT$ in the triangulation we
  choose an affine linear function $\chi_\tau: A\rightarrow \R$ with the
  following property: 
 $\chi_\tau(\omega)\geq\lambda(\omega)$ with equality holding exactly
 for $\omega\in\tau\cup\orr$. Note that for $\sigma$, an
 $(N-1)$-dimensional simplex, $\chi_\sigma$ is uniquely
 determined by $\psi_\lambda\arrowvert_{C_\sigma}$ and the inequality
 condition is satisfied because $\psi_\lambda$ is a strictly convex
 function. By the same reason we can satisfy the inequality for the
 simplices of smaller dimension. For every  $(N-1)$-dimensional
  simplex $\sigma$ we define the function 
$$\chi_\sigma(s,\gamma,\omega):=-\log(\sum_{\tau\in\dT}\rho^0_\tau
 e^{-\gamma\cdot\chi_\tau}),$$
 first for $\omega\in\sigma\cup\orr$ and then extend by linearity to all
 $\omega\in A$.   

The function $\chi(s,\gamma,\omega)$ is constructed
  by gluing the functions $\chi_\sigma(s,\gamma,\omega)$ together in
  the following way.
The collection of maximal dimensional simplices $\sigma\in\dT$
provides a triangulation of $\dD$. We take small open neighborhoods of
each $\sigma$ to get an open cover $\{\widetilde{\sigma}\}$ of $\dD$
and choose a partition of unity $\{\alpha_\sigma\}$ subordinate to
it. We require the open enlargements of $\sigma$'s to be small enough, so
that $\widetilde{\sigma}\subset\bigcup_{\tau\subset\sigma}W_\tau$. In
particular, $\rho^0_\tau\arrowvert_{\widetilde{\sigma}}\equiv 0$ unless
$\tau\subset\sigma$. Now we can define the desired function
$$\chi(s,\gamma,\omega):=\sum_{\sigma\in\dT}\alpha_\sigma(s)
\chi_\sigma(s,\gamma,\omega),$$ 
which is smooth with respect to
   $(s,\gamma)$ and affine linear with respect to $\omega$ by
   construction. It is also clear that
  $$\chi_\sigma(s,0,\omega)=-\log(\sum_{\tau\in\dT}\rho^0_\tau)=0$$ and
   $$\chi_\sigma(s,\gamma,\orr)=-\log(\sum_{\tau\in\dT}\rho^0_\tau
   e^{-\gamma\cdot\lambda(0)})=\gamma\cdot\lambda(0)$$ for all
   $\sigma$. Hence
   $\chi(s,0,\omega)\equiv 0$, and  $\chi(s,\gamma,\orr)\equiv
   \gamma\cdot\lambda(0)$.  

The last thing to check is the behavior of $\chi(s,\gamma,\omega)$ as
$\gamma\rightarrow\infty$. Fix a 
point $s\in W_{\tau_0}$, and first consider vertices $\omega\in\tau_0$.
The partition functions $\alpha_\sigma(s)=0$ unless
$\tau_0\subset\sigma$, hence $\omega\in\sigma$ and 
$\chi(s,\gamma,\omega)= -\log(\sum_{\tau\in\dT}\rho^0_\tau
 e^{-\gamma\cdot\chi_\tau(\omega)})$.
$$\lim_{\gamma\rightarrow\infty} e^{-\chi(s,\gamma,\omega)}
 e^{\gamma\cdot\lambda(\omega)}=
 \lim_{\gamma\rightarrow\infty}\sum_{\tau\in\dT}\rho^0_\tau
 e^{\gamma\cdot(\lambda(\omega)-\chi_\tau(\omega))}=
 \sum_{\tau\ni\omega}\rho^0_\tau=\rho_\omega,$$
because $\lambda(\omega)-\chi_\tau(\omega)\leq 0$ with equality
 holding exactly for $\omega\in\tau\cup\orr$.  

If $\omega\notin\tau_0$, then to show that
$\lim_{\gamma\rightarrow\infty} e^{-\chi(s,\gamma,\omega)}
e^{\gamma\cdot\lambda(\omega)}=\rho_\omega(s)=0$ we look at the
asymptotics of 
the affine functions $\chi_\sigma(s,\gamma,\omega)$ as
$\gamma\rightarrow\infty$. Let   $Q$ be
 the collection of simplices $\tau$ with the minimal value of
 $\chi_\tau(\omega)$ among those with nonzero $\rho_\tau^0(s)$. Note that if
 $\rho_\tau^0(s)\neq 0$, then $\tau\subseteq\tau_0$, and hence
 $\chi_\tau(\omega)>\lambda(\omega)$ as $\omega\notin\tau$.
So we see that
$$\chi_\sigma(s,\gamma,\omega)\sim -\log(\sum_{\tau\in Q} \rho^0_\tau) +
 \gamma\cdot\chi_\tau(\omega).$$
Combining these together for all $\sigma$'s we get 
$$\chi(s,\gamma,\omega)\sim -\log(\sum_{\tau\in Q} \rho^0_\tau) +
 \gamma\cdot\chi_\tau(\omega),$$
and hence 
$$ e^{-\chi(s,\gamma,\omega)}e^{\gamma\cdot\lambda(\omega)}
\sim\sum_{\tau\in Q}\rho^0_\tau 
e^{\gamma\cdot(\lambda(\omega)-\chi_\tau(\omega))}\rightarrow 0$$
as $ \gamma\rightarrow\infty$. This completes the proof.
  \end{proof}
 \end{lemma}
 
 \begin{remark}
   In the proof of the above lemma the crucial fact we used was that
   the characteristic function $\psi_\lambda$ is strictly convex.  
 \end{remark}
Just as we did for $\rho_\omega$, we will use the same notation for both
$\chi(\gamma,\omega)$ and its pull back to $X_\Delta^0$ via the moment
map. Now we have all the tools to prove the main theorem. Let $R$ be the
positive real number as in proposition 3.1. Denote by $H_R$ the
one-(real) parameter family of the auxiliary hypersurfaces
$\{H_t,|t|=R\}$. Let $F_{\Gamma_0 R}$ be the one-(complex) parameter
family of the original Calabi-Yau hypersurfaces $\{F_ t,|t|\geq
\Gamma_0 R\}$.
\begin{thm}
   There is a positive real number $\Gamma_0$, such that there exists
  a diffeomorphism between the families $\{H_R\}\times 
  (\Gamma_0,+\infty)$ and $F_{\Gamma_0 R}$, which specializes to a
  diffeomorphism between 
  the hypersurfaces $(H_t,\Gamma)$ and $F_{\Gamma t}$.
  \begin{proof}
    First we define a hypersurface  $H^\varepsilon_t\subset\X^0$ by
    the equation 
    $$ t^{\lambda(0)} x^{\orr}-\sum_{\omega\in A\cap\dD}(\rho_\omega
    +\varepsilon_\omega) t^{\lambda(\omega)}x^\omega =0,\quad \text{where
    } \varepsilon_\omega :=
    e^{\gamma\lambda(\omega)-\chi(\gamma,\omega)}-\rho_\omega. $$
    According to lemma 3.4, all $\varepsilon_\omega$ uniformly vanish as
    $\gamma\rightarrow\infty$, that is for $\gamma\geq\gamma_0$
    $H^\varepsilon_t$ is
    indeed a small deformation of $H_t$ and hence diffeomorphic to
    it. Using the substitution  
    ${x'}^\omega:= x^\omega e^{-\chi(\gamma,\omega)}$ we
    get
$$  e^{\gamma\lambda(0)} t^{\lambda(0)}{x'}^{\orr}-\sum_{\omega\in
    A\cap\dD}e^{\gamma\lambda(\omega)} t^{\lambda(\omega)}{x'}^\omega
    =0, $$ which is exactly the equation of the hypersurface
    $F_{(\Gamma\cdot t)}$ for $\Gamma=e^\gamma$. A priori this
    substitution defines a map
    $\X^0\rightarrow\PP^{|\Delta(\Z)|-1}$. But because $\chi(\omega)$
    is an affine function, the image of this map, in fact, lies in
    $\X$. Hence the above equation indeed defines a hypersurface in
    $\X$.

    Notice that this construction works for the entire families,
    because the deformation diffeomorphisms clearly form a trivial
    system, and the substitutions depend only on the absolute values of
    the parameters of the families.
  \end{proof}
\end{thm}

 According to [GKZ], Ch. 10, all hypersurfaces which lie inside a
 translated  cone $C(T)+b$, where $b$ is some vector in the interior
 of $C(T)$, are smooth and hence diffeomorphic to each other.
 Combining the above theorem with proposition 3.2, we get the main
 result of the paper.

\begin{cor}
  A Calabi-Yau hypersurface in $\X$, which is sufficiently far away 
 from the walls of the secondary fan to $\Delta$ and
 sufficiently close to the large complex structure, admits a fibration
 over a sphere $S^{N-1}$ with generic fibers $(N-1)$-dimensional
  tori.
\end{cor}

 \begin{remark}
   It should be possible to remove the smoothness
   requirement. In this case one has to be more careful with
   deforming the equation of a non smooth hypersurface. For a
   $\Delta$-regular hypersurface in the translated cone $C(T)+b$ all
   singularities come from the singularities of $\X$. There is a
   natural stratification of $\X$ by $\mu^{-1}(\Theta)$, as $\Theta$ runs over
   open parts of the faces of $\Delta$, which induces a stratification
   of $F_t$. All diffeomorphisms should
   then be understood in this stratified sense (see, e.g. [GM].)  
 \end{remark}

The ultimate goal would be, of course, to construct a special
Lagrangian fibration. Our construction, unfortunately, leaves this
problem open. But there are some features which may be worth
mentioning. For instance, our fibration is quite special in the
following sense. It 
tends to concentrate the singularities of the fibers into a smaller
number of fibers with worse degenerations. As an example let us consider
the family of quartic K3 surfaces in $\C\PP^3$ given by the equations 
$$t\cdot x^{\orr}+ \sum_{\Omega\text{ vertex of }\Delta} x^\Omega +
O(t^{-1})=0.$$ 
A generic fibration is expected to have 24 degenerate fibers, and each one of
them is homeomorphic to the standard $I_1$ degenerate elliptic curve. In
our fibration the terms $O(t^{-1})$ don't matter and we get just 6
singular fibers of type $I_4$. 

There is a local special Lagrangian structure on the algebraic pieces
of the auxiliary hypersurface. However for this we should have defined
the cutoff functions $\rho_\omega$ slightly 
differently (we didn't do so in the first place because it would have
spoiled the uniformness of the definition). Namely, let
$\rho_\omega:=\sum_{\tau}\rho^0_\tau$, where $\tau$ runs over
the simplices in $\dT$ containing $\omega$, but of dimension
at most $(N-2)$. This reduces the support 
of $\rho_\omega$ to a neighborhood of the $(N-2)$-skeleton of
$\dD$. Then for a maximal dimensional simplex $\sigma\subset\dT$ the
equation of $H_t$ in $\U$ would be
just $x^{\orr}=0$, which defines some open subset of the corresponding
to $\sigma$ toric divisor. Notice that $\mu_t(T_s)$ is just one point
in $\Delta$ (which in this case lies on the boundary $\dD$), so with
respect to the standard symplectic form on the 
toric variety $\X$, $T_s$ is clearly Lagrangian. Moreover, if we define
the top holomorphic form in $\widetilde{U_\sigma}$ according to the
equation (every hypersurface is locally Calabi-Yau), then it restricts
to a volume form on each fiber. It is easy to check that the same is true
in $\widetilde{U_\omega}$ for any vertex in the triangulation, where the
local equation is $y_1...y_N=const$. The fibration is given by fixing
$|y_i|, 1\leq i\leq N$, which is clearly Lagrangian. A top holomorphic
form can be written as
$\frac{dy_1}{y_1}\wedge...\wedge\frac{dy_{N-1}}{y_{N-1}}$, which
  restricts to a volume form on $T_s$. 

 So that in
the right $W$-decomposition limit every $T_s$ becomes Lagrangian with
respect to the deformed symplectic structure except for $s$ in
the singular locus. Unfortunately we cannot say the
same thing about special Lagrangian property. Although the local
holomorphic forms on the auxiliary hypersurface do give the volume
forms when restricted to the 
fibers, it is not at all clear what are their pull-backs to the
original hypersurface and how to patch them together in the
transition regions.

\section{Monodromy}

In this section we want to show an application of the
constructed fibration to the monodromy calculations.
Gross has made a conjecture [G] about the monodromy transformation in a
family of Calabi-Yau manifolds ${\mathcal X}\rightarrow S$. Let
$X={\mathcal X}_t\stackrel{f}{\rightarrow} B$, $t\in S$ be a torus
fibration with a section $\delta_0$. Then $X^\sharp$, the complement of the
critical locus of $f$, has a structure of a fiber space of abelian groups
with the zero section $\delta_0$. Given another section $\delta$ one
obtains a diffeomorphism $T_\delta:X^\sharp\rightarrow X^\sharp$ given by
$x\mapsto x+\delta(f(x))$, which extends to a diffeomorphism of the
entire $X$. Given a degeneration divisor in $S$ passing
through the large complex structure 
point, and a loop around this divisor we can consider the monodromy
transformation on cohomology. The conjecture says that this monodromy
is induced by $T_\delta$ for some section $\delta$. 

We are going to construct such a section for our family of
hypersurfaces $F_t$. Because theorem 3.5 establishes the
diffeomorphism between the entire families, the monodromy question is
identical for the auxiliary family $H_R$. Without loss of generality
we may assume that the base point in the family is given by a hypersurface
$H_{t_0}$ with $t_0=R$, a real positive number. The monodromy loop is
parameterized by $t=t_0 e^{2\pi i \gamma},\ 0\leq\gamma\leq1$. But
first of all we need a zero section.

\begin{lemma}
  The fibration $h_{t_0}:H_{t_0}\rightarrow\dD$ has a section $\delta_0$,
  which misses all singular points of the fibers.
  \begin{proof}
    The section is given by the set of all real positive points of
    $H_{t_0}$. Let $s\in W_\tau$ be a $\tau$-point. We just have to
    show that the interval $I_s$ of positive real points in $\X^0$ has
    a unique solution to the equation
$$    t^{\lambda(0)}x^{\orr}=\sum_{\omega\in 
  \tau}\rho_\omega  t^{\lambda(\omega)}x^\omega.$$
   But with the identification of the real positive points of $\X^0$
  with $\Delta^0$, the real positive points satisfying this equation
  form a hypersurface separating 
  $\orr$ from those $\omega$ for which $\rho_\omega\neq 0$
  (cf. [GKZ], Ch.11). In particular, it has a unique point of intersection
  with the line $I_s$, which, moreover, lies in the interior of
  $\Delta^0$. But all the singular points of the fiber $T_s$ are
  mapped to the boundary of $\Delta$. Thus we get the desired section
  $\delta_0: \dD\rightarrow H_{t_0}$.
  \end{proof}
\end{lemma}

To construct the other section $\delta$ we consider a Delzant type polytope
$\DD\in M^*$ defined by the inequalities:
$$\langle m ,\omega-\orr\rangle\geq\gamma\cdot(\lambda(\omega)-\lambda(0)),$$
where $\omega$ runs over the vertices in the triangulation $T$.
This is a convex polytope with non empty interior for $\gamma>0$,
because of strict convexity of the characteristic function
$\psi_\lambda$. By the same reason it is combinatorially dual to
$(\Delta,T)$, namely  to each $k$-dimensional simplex $\tau$ in $\dD$,
there corresponds an $(N-1-k)$-dimensional face $\tau^\vee$ in
$\partial\Delta^\vee_\lambda$ with the reverse incidence relation.

 The bijective correspondence between the centers of the dual
 pairs, $\tau$ and $\tau^\vee$, gives rise to a simplicial map 
$\nu_\gamma: Bar(\dT)\rightarrow Bar(\dDD)$ between the first barycentric
subdivisions. Considering the $W$-decompositions, which are dual to
the barycentric ones, we get a homeomorphism $\nu_\gamma:\dD
\rightarrow\dDD$ satisfying
$\nu_\gamma(W_\tau)=W_{\tau^\vee}$, where $\{W_{\tau^\vee}\}$
provide a CW-decomposition of $\dDD$. Because each $W_{\tau^\vee}$
contains the center of the simplex $\tau^\vee$, there is a map
$\nu'_\gamma:\dD \rightarrow\dDD$, homotopic to $\nu_\gamma$,
with the property that $\nu'_\gamma(W_\tau)\subset\tau^\vee$.

We will use the same notation for both $\nu'_\gamma$ and its pull back to
$\X^0$. Let us define a diffeomorphism $D_\gamma:\X^0\rightarrow\X^0$ by 
$$ x^\omega\mapsto x^\omega e^{2\pi i\<\nu'_\gamma(x),\omega-\orr\>}.$$
This is well defined as $\nu'_\gamma(x)$ is a linear functional with respect
to $\omega$. In fact, this diffeomorphism is equivariant with respect
to the toric action, i.e. $\mu_t\circ D_\gamma =\mu_t$. For
convenience we will drop the index $\gamma$ in all notations whenever
$\gamma=1$. The desired section $\delta:\dD\rightarrow H_{t_0}$ is
given by applying the diffeomorphism $D$ to the zero section. Thus we
define $\delta:=D\circ\delta_0$.
\begin{thm}
  The section $\delta:\dD\rightarrow H_{t_0}$ is well defined and induces
  the monodromy transformation on $H_{t_0}$.
  \begin{proof}
    The correctness of the definition follows easily from the
    following observation. For  $s\in W_\tau$ notice that
    $\langle\nu'_\gamma(s),\omega-\orr\rangle =
    \gamma(\lambda(\omega)-\lambda(0))$ for all $\omega\in\tau$. This
    means that for $\gamma=1$ the diffeomorphism $D$ has no effect on
    any monomial $x^\omega$ for $\omega\in\tau$, as it gets 
  multiplied by the factor of $e^{2\pi i(\lambda(\omega)-\lambda(0))}$. So that
    the equation of $H_{t_0}$ in $\W$  
$$ t_0^{\lambda(0)}x^{\orr}-\sum_{\omega\in
  \tau}\rho_\omega  t_0^{\lambda(\omega)} x^\omega =0$$
  is still satisfied for $\delta=D(\delta_0)$. Notice that the
  diffeomorphism $D$ respects the fibration
  $h_{t_0}:H_{t_0}\rightarrow\dD$, hence the action in a fiber $T_s$
  is just the translation by $\delta(s)-\delta_0(s)$.  This action is
  also well defined on singular fibers. Indeed,  a singular fiber is
  itself a fibration $p_s:T_s\rightarrow T^{\dim\tau}$, according to the
  proposition 3.2. And the action on $T_s$ translates points along
  the fibers of $p_s$, which generically are $(N-1-\dim\tau)$-dimensional
  tori. 

  To see that the diffeomorphism $D$ induces the monodromy
 transformation we notice that $D_\gamma$ provides a  diffeomorphism
 of $H_{t_0}$ with $H_{t_0 e^{2\pi i \gamma}}$. Indeed, in $\W$ the
 defining equation of $H_{t_0}$ translates into
$$ t_0^{\lambda(0)}x^{\orr}-\sum_{\omega\in
  \tau}\rho_\omega  t_0^{\lambda(\omega)} x^\omega e^{2\pi
  i\gamma(\lambda(\omega)-\lambda(0))} =0,$$ which is exactly the defining
  equation for $H_{t_0 e^{2\pi i \gamma}}$. As $\gamma$ runs from 0 to
  1 the family of the diffeomorphisms $D_\gamma$ provides the
  monodromy  along the loop $t_0 e^{2\pi i\gamma}$.
 This completes the proof.
  \end{proof} 
\end{thm}

\section {Dual fibrations and mirror symmetry}

This section is rather speculative in character but it is 
impossible to overlook a connection of our construction
with the mirror symmetry. A triple $(\Delta,T,\lambda)$ defines a
family of complex structures on a Calabi-Yau hypersurface. For
simplicity we assume that the subset $A\subset\Z^N$ coincides with the
set of vertices of the triangulation $T$. On the mirror side we want to
get a family of K\"ahler structures on some other Calabi-Yau. This
family is provided by the monomial-divisor map [AGM]. To construct it
we consider the polytopes $\DD$ defined in the previous section. 

Let $N(\DD)$ be the normal fan to $\DD$.  $N(\DD)$ is a rational convex
polyhedral fan and the corresponding toric
variety $X_{\DD}$ is a blow-up of the variety $X_{\Delta^D}$
for the polar polytope $\Delta^D$. The one-dimensional
cones in $N(\DD)$ are in one-to-one correspondence with the
vertices  of the triangulation $\dT$. The exceptional
divisors are labeled by those of them which are not the vertices of
$\Delta$ and  $N(\DD)$ is a simplicial cone subdivision of
$N(\Delta^D)$ by means of the triangulation $T$. The vector $\lambda$
lies in the interior of the K\"ahler
cone of $\XX$ and defines the K\"ahler class (in the orbifold sense)
by the linear combination of the toric divisors $[\omega^\vee]$
corresponding to the facets $\omega^\vee$ in $\XX$: 
$$[\kappa_\gamma]:=-\sum_{\omega\in\dT}\gamma(\lambda(\omega)
-\lambda(0))[\omega^\vee].$$ 
We will consider the family $\{\XX\}$, where $\gamma$
runs over positive real numbers. The symplectic form $\kappa_\gamma$ defines
the moment map
$\mu_\gamma^\vee:\XX\rightarrow\DD$ (cf. [Gu]). Now we choose a
regular anti-canonical hypersurface $Z^D$ in $X_{\Delta^D}$ with large
complex structure (e.g., with
large central coefficient). Denote by $Z\subset\XX$ its proper
transform induced by the blow-up $\XX\rightarrow X_{\Delta^D}$. $Z$
is a Calabi-Yau hypersurface (cf. [B]) endowed with a K\"ahler
structure by restriction from $\XX$ in the orbifold sense. This is the 
mirror family. We will let $\gamma=1$ for the future consideration as the
behaviour of the family changes by a simple rescaling for other
$\gamma$.

To study the geometry of $Z$ we will again use the moment map
$\mu^\vee: X_{\Delta^\vee}\rightarrow\Delta^\vee$. At this point we
will make an assumption that
$Z$ possesses a torus fibration analogous to that of a smooth
hypersurface. All singularities of $Z$ are mapped by $\mu^\vee$ to the
$(N-2)$-skeleton of $\Delta^\vee$. So that a generic fiber is
still a smooth $T^{N-1}$. But degenerations in singular tori may give
rise to the singularities in the total space.

First, let us introduce some notations. We will use $\tau_\Z$ also to
denote the subgroup of $\Z^N$ modeled on the affine
sublattice $\tau_\Z$. In other words, $\tau_\R:=\tau_\Z\otimes\R$ is the
$k$-dimensional vector subspace parallel to $\tau$ and passing through
$\orr$. Denote by
$\tau_\Z^*$ the quotient of $(\Z^N)^*$ dual to $\tau_\Z\subset\Z^N$,
and let $\tau_\R^*:=\tau_\Z^*\otimes\R$ be the corresponding quotient
of $M^*$. Now we consider an explicit parameterization of nonsingular
fibers in the original family. Let $s\in W_\tau$,
for $\tau\subset\dT$, be a point in $\dD$. According to
proposition 3.2, a fiber $T_s$ is a fibration itself $p_s:T_s\rightarrow T^k$,
and $T^k$ is naturally isomorphic to the torus
$\tau^*_\R/\tau^*_\Z$. Choosing a point in $T^k$ determines the phases
not only of $x^\omega,\ \omega\in\tau$, but also of
$x^{\orr}$. Considering the fact that $\tau^\vee_\R$ is defined by the
equations  $\langle u,\omega-\orr\rangle =0$, $\omega\in\tau$, we
conclude that an $(N-k-1)$-dimensional fiber $T^{N-k-1}$ can be
identified with $\tau^\vee_\R/\tau^\vee_\Z$. Hence the fiber $T_s$ is
isomorphic to the torus $\tau^*_\R/\tau^*_\Z\oplus
\tau^\vee_\R/\tau^\vee_\Z$, though the splitting into the direct sum is
not natural. The dual fiber $T_{s^\vee}$, where $s^\vee=\nu(s)$ is a
point in $W_{\tau^\vee}\subset\dD^\vee$, will be isomorphic to
$\tau_\R/\tau_\Z\oplus (\tau^\vee_\R)^*/(\tau^\vee_\Z)^*$. 

To conclude the picture we need to take into consideration the
singular tori. Remember that the discriminant locus $D(H_t)$ is
homotopy equivalent (and, in fact, can be made arbitrarily close in
the appropriate $W$-decomposition limit)
 to $Sk_T^{(N-3)}$, the $(N-3)$-skeleton of the subdivision dual to
the triangulation of $\partial\dD$. $Sk_T^{(N-3)}$ is a simplicial
complex consisting of the
simplices $(O(\tau_{i_1}),...,O(\tau_{i_k}))$, with vertices
$O(\tau_{i_j})$, the centers of $\tau_{i_j}$, and 
$\tau_{i_1}\subset...\subset\tau_{i_k}$ 
running over all nested chains of simplices in $\partial\dD$ of positive
dimension.

On the mirror side the discriminant locus is again homotopy equivalent
to a simplicial complex  $(Sk_T^{(N-3)})^\vee$ with the vertices
$O(\tau^\vee)$ and the simplices labeled by the nested chains of
$\tau^\vee$'s. However the simplices
$\tau^\vee\subset\partial\dD^\vee$ which have appeared as a 
result of the blow up $X_{\Delta^\vee}\rightarrow X_{\Delta^D}$ do not
contain any points in the image of the moment map $\mu^\vee(Z)$, hence
should be excluded from the discriminant locus. They correspond
exactly to those simplices $\tau\in \dT$, for which the minimal face
$\Theta_\tau$ is a facet of $\Delta$, i.e. to those which are not in
$\partial\dD$. 

The simplicial map $\nu:\dD\rightarrow\dD^\vee$ provides a one-to-one
correspondence between the points $O(\tau)$ and  $O(\tau^\vee)$, and
thus establishes the simplicial isomorphism between  $Sk_T^{(N-3)}$
and  $(Sk_T^{(N-3)})^\vee$. So that in the right limit we get the
identification between the two discriminant loci. This suggests to
consider the corresponding singular fibers $T_s$ and $T_{\nu(s)}$ to
be dual to each other.

\end{document}